\documentclass[12pt]{amsart}
\usepackage{amssymb,amsmath,amsthm}
\usepackage{amsfonts}
\usepackage[dvips]{graphics}
\usepackage{epsfig}
\usepackage{graphicx}
\usepackage{psfrag}
\usepackage{amsmath}
\usepackage{amssymb,subfigure}

\def\R {\mathbb{R}}

\def \and{\quad\text{and}\quad}

\def\ds{\displaystyle}

\def\R{\mathbb R}

\def\epsilon{\varepsilon}
\def\ds{\displaystyle}

\newcommand{\be}{\begin{equation}}
\newcommand{\ee}{\end{equation}}
\newcommand{\baa}{\begin{array}}
\newcommand{\eaa}{\end{array}}
\newcommand{\ba}{\begin{eqnarray}}
\newcommand{\ea}{\end{eqnarray}}

\numberwithin{equation}{section}

\newtheorem{a2theorem}{Theorem}[section]

\newtheorem{a2remark}{Remark}[section]

\def \no#1#2#3 {{\bf #1} (#3), #2.}
\def \eds#1#2#3 {#1, #2, #3.}

\theoremstyle{definition}

\numberwithin{equation}{section}

\title[Population dynamics under the influence of external
perturbations] {Models of population dynamics under the influence of
external perturbations: mathematical results}

\author[M.D. Chekroun, and L. Roques]
{Micka\"el D. Chekroun and Lionel Roques}

\address{\'Ecole Normale Sup\'erieure - CERES-ERTI
\newline\indent
75005 Paris, France} \email{chekro@lmd.ens.fr {\rm (M. D. \
Chekroun)}}

\address{Unit\'e Biostatistique et Processus Spatiaux, INRA,
Domaine St Paul - Site Agroparc, 84914 Avignon Cedex 9, France}
\email{lionel.roques@avignon.inra.fr {\rm (L. Roques)}}

\thanks{The present manuscript has been published as:
\textsc{M.D. Chekroun, and L. Roques}, Models of population dynamics
under the influence of external perturbations: mathematical results,
{\it C. R. Acad. Sci. Paris}, Ser. I, 343: 307--310.}



\begin{document}
\maketitle

\begin{abstract}
In this Note, we describe the stationary equilibria and the
asymptotic behaviour of an heterogeneous logistic reaction-diffusion
equation under the influence of autonomous or time-periodic forcing
terms. We show that the study of the asymptotic behaviour in the
time-periodic forcing case can be reduced to the autonomous one, the
last one being described in function of the ``size" of the external
perturbation. Our results can be interpreted in terms of maximal
sustainable yields from populations. We briefly discuss this last
aspect through a numerical computation.



\vspace{2ex}
\noindent{\bf R\'esum\'e:} {\bf Analyse de mod\`eles de dynamique de
populations sous l'influence de perturbations externes. } Cette Note
a pour objet l'\'etude des \'etats stationnaires et du comportement
asymptotique d'\'equations de r\'eaction-diffusion avec coefficients
h\'et\'erog\`enes en espace, auxquelles nous ajoutons un terme de
perturbation stationnaire ou p\'eriodique en temps. Nos r\'esultats
peuvent s'interpreter en termes de r\'ecolte maximale supportable
par une population. Nous soulignons cet aspect \`a l'aide d'un
calcul num\'erique.

\end{abstract}



\section{Introduction}\label{sec1}

The purpose of this Note is to study the following model:
\be\label{eq_evo} u_t= \nabla \cdot(A(x) \nabla u)+u(\mu(x)-\nu(x)
u)- f(\omega t,x)\rho_{\varepsilon}(u), \ (t,x) \ \in \
\R_+\times\Omega. \ee The reaction-diffusion models of the type
$u_t= \nabla \cdot(A(x) \nabla u)+u(\mu(x)-\nu(x) u)$ correspond to
the natural extension of the classical Fisher model \cite{fi}. They
were first introduced by Shigesada et al. \cite{skt} for population
dynamics. Our aim is to understand the asymptotic behaviour of the
solutions of such models, when we add a time-periodic forcing term
$f(\omega t,x)$. With such additional term, this can interpreted as
an \emph{harvesting model} with seasonal harvesting. In real-life
context this
 perturbation term can arise when a quota is set on the harvesters.

We make the following assumptions on the coefficients: the diffusion
matrix $A(x)$ is assumed to be of class $C^{1,\alpha}$ (with
$\alpha>0$) and uniformly elliptic; i.e. there exists $\tau>0$ such
that $A(x)\ge\tau I_N$ for all $x\in\Omega$. The functions $\mu$ and
$\nu$ belong to $L^{\infty}(\Omega)$. Moreover, we assume that there
exist $\underline{\nu}$ and $\overline{\nu}$ such that
$0<\underline{\nu}<\nu(x)<\overline{\nu}$ for all $x$ in $\Omega$.
The function $f$ is 1-periodic in the first variable and belongs to
$C^0(\R\times\Omega)$, and the function $\rho_{\varepsilon}$ defines
a ``regularized threshold'': it is a $C^1(\R)$ nondecreasing
function such that $\rho_{\varepsilon}(s)=0$ for all $s\leq 0$ and
$\rho_{\varepsilon}(s)=1$ for all $s\geq \varepsilon$. This
threshold guarantees the non-negativity of the solutions of
(\ref{eq_evo}).

Two kinds of domains $\Omega$ are considered: either $\Omega=\R^N$
or $\Omega$ is a smooth bounded domain of $\R^N$. We qualify the
first case, $\Omega=\R^N$,  as the {\em sp-case} and the second one
as the {\em bounded case}. Indeed, in the sp-case, we assume that
$A(x)$, $\mu(x)$, $\nu(x)$ and $f(s,x)$ depend on the variables
$x=(x_1,\cdots,x_N)$ in a space-periodic fashion (i.e. for
$L_1,\ldots,L_N$ fixed positive numbers, a function $g$ is said to
be sp-periodic if $g(x+k)=g(x)$ for all $x\in\mathbb{R}^N$ and $k\in
L_1\mathbb{Z}\times\cdots\times L_N\mathbb{Z}$). In the bounded
case, throughout this paper, we assume that we have Neumann boundary
conditions on $\partial \Omega$.

\section{The case of autonomous forcing}\label{sec2}
All the results of this section remain true either in the {\bf
sp-periodic} or {\bf bounded cases}. The proofs are detailed in
\cite{rc}.

We consider the equation (\ref{eq_evo}) with $f(wt,x)=\delta h(x)$,
i.e. \be\label{eq_evo2} u_t= \nabla \cdot(A(x) \nabla
u)+u(\mu(x)-\nu(x) u)-\delta h(x)\rho_{\varepsilon}(u), \ (t,x) \
\in \ \R_+\times\Omega, \ee where $h$ is a continuous function such
that there exist $\alpha,\beta>0$ with $\alpha< h(x)<\beta \hbox{
for all }x \in \Omega$, and which is sp-periodic in the sp-case.

Let $\lambda_{1}$ be defined as the unique real number such that
there exists a function $\phi>0$ which satisfies
\be\label{eqlambda1} -\nabla \cdot(A(x) \nabla
\phi)-\mu(x)\phi=\lambda_{1}\phi\hbox{ in }\Omega, \phi>0 \hbox{ and
 }\|\phi\|_{\infty}=1, \ee with either periodic or Neumann boundary
conditions, depending on $\Omega$, as mentioned above. The function
$\phi$ is uniquely defined by (\ref{eqlambda1}) (the existence and
uniqueness of $\lambda_1$ and $\phi$ follow from the standard
Krein-Rutman theory).

\vspace{1ex}
\begin{a2remark} \label{rem1}{\rm Note that if we assume that $\lambda_1<0$ and $\delta=0$, then, given any
continuous and bounded function $u_0$, the solution $u(t,x)$ of
(\ref{eq_evo2}) with initial data $u_0$ converges to a function $p$
which is the unique bounded and positive  solution of $\nabla
\cdot(A(x) \nabla p)+p(\mu(x)-\nu(x) p)=0$, $x \in\Omega$. These
convergence,  as well as existence and uniqueness results are proved
in \cite{bhr1}.}
\end{a2remark}
\vspace{1ex}


We first describe the steady states of (\ref{eq_evo2}) without
``regularized threshold'': \be\label{eq_sta} \nabla \cdot(A(x)
\nabla p_{\delta})+p_{\delta}(\mu(x)-\nu(x) p_{\delta})-\delta
h(x)=0, \ x \ \in \ \Omega.\ee Using a Leray-Schauder degree
argument, together with the uniqueness of the solution $p$ defined
in the above remark, we prove the following\hfill\break

\begin{a2theorem}
\label{CRASth1} There exists $\delta^*>0$ such that for all $\delta$
s. t. $0<\delta<\delta^*$, {\rm(\ref{eq_sta})} admits  two distinct
positive solutions, $p^1_\delta$ and $p^2_\delta$. Moreover,
$p^1_\delta\to 0$ and $p^2_\delta\to p$ uniformly in $\Omega$ as
$\delta\to 0$.
\end{a2theorem}

\vspace{1ex}

Let us set $\ds{\underline{\phi}:=\min_{x\in \Omega}\phi(x)},$  $
\delta_1:=\ds{\frac{\lambda_1^2 \underline{\phi}}{\beta
\overline{\nu}(1+\underline{\phi})^2}}$ and
 $\delta_2:=\ds{\frac{\lambda_1^2 }{4 \alpha \underline{\nu}}}.$
 Then we have the following theorem:


\vspace{1ex}

\begin{a2theorem}\label{CRASth2}
\par - (i) If $\lambda_1<0$ and $\delta\leq \delta_1$, then there
exists a positive bounded solution $p_\delta$ of
{\rm{(\ref{eq_sta})}} such that $p_\delta \geq \ds{-\frac{\lambda_1
\phi}{\overline{\nu}(1+\underline{\phi})}} $ (in particular $\max
p_\delta \geq \frac{-\lambda_1}{2\overline{\nu}}$).
\par (ii) If $\lambda_1<0$ and $\delta
> \delta_2$, or if $\lambda_1\geq 0$,  there is no positive bounded
solution of {\rm{(\ref{eq_sta})}}.
\end{a2theorem}

\vspace{1ex}

\begin{figure} \centering
\subfigure[]{%
\label{fig+2} 
\includegraphics*[width=4.1cm]{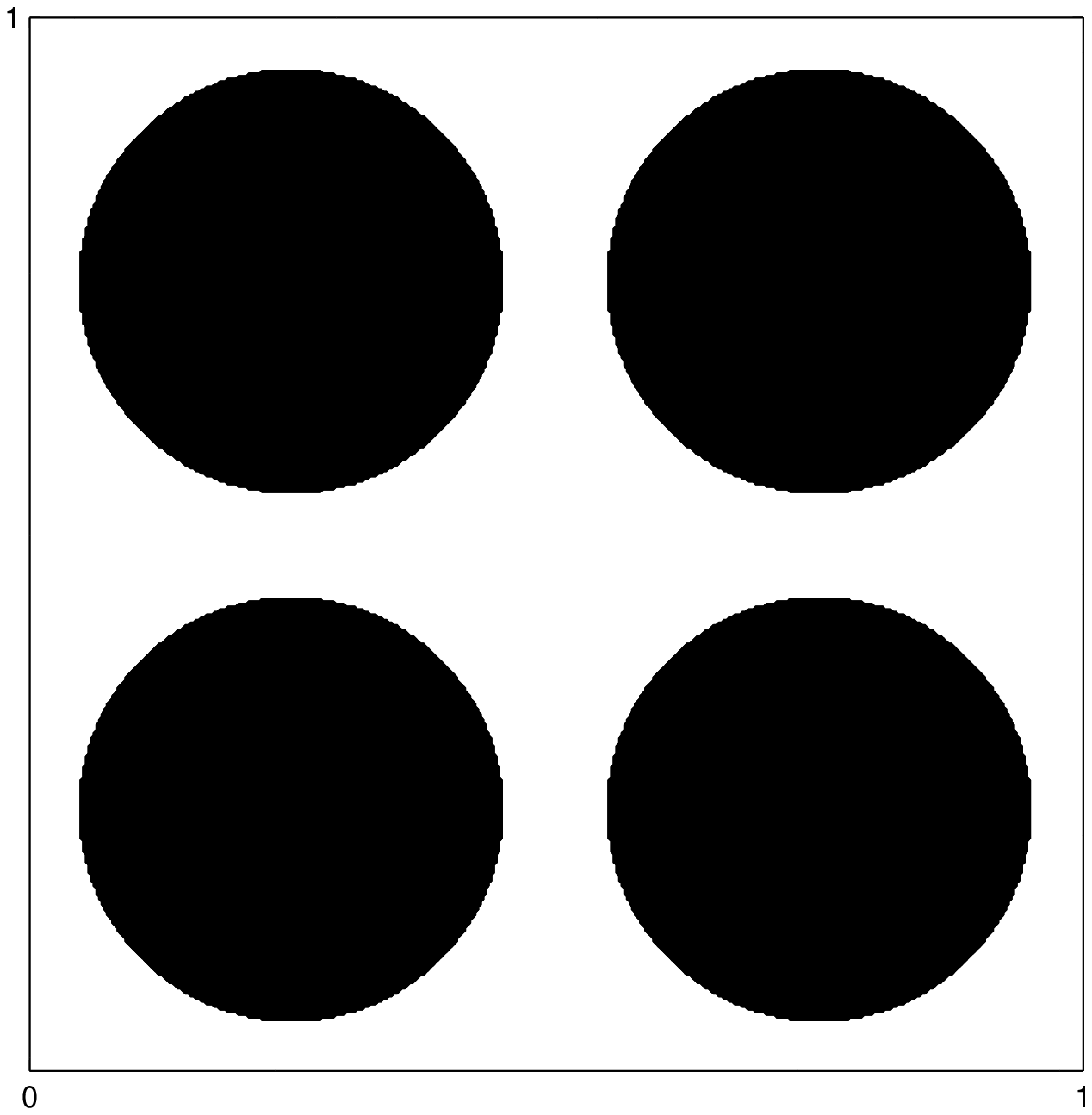}}
\subfigure[]{%
\label{fig+10} 
\includegraphics*[width=4.1cm]{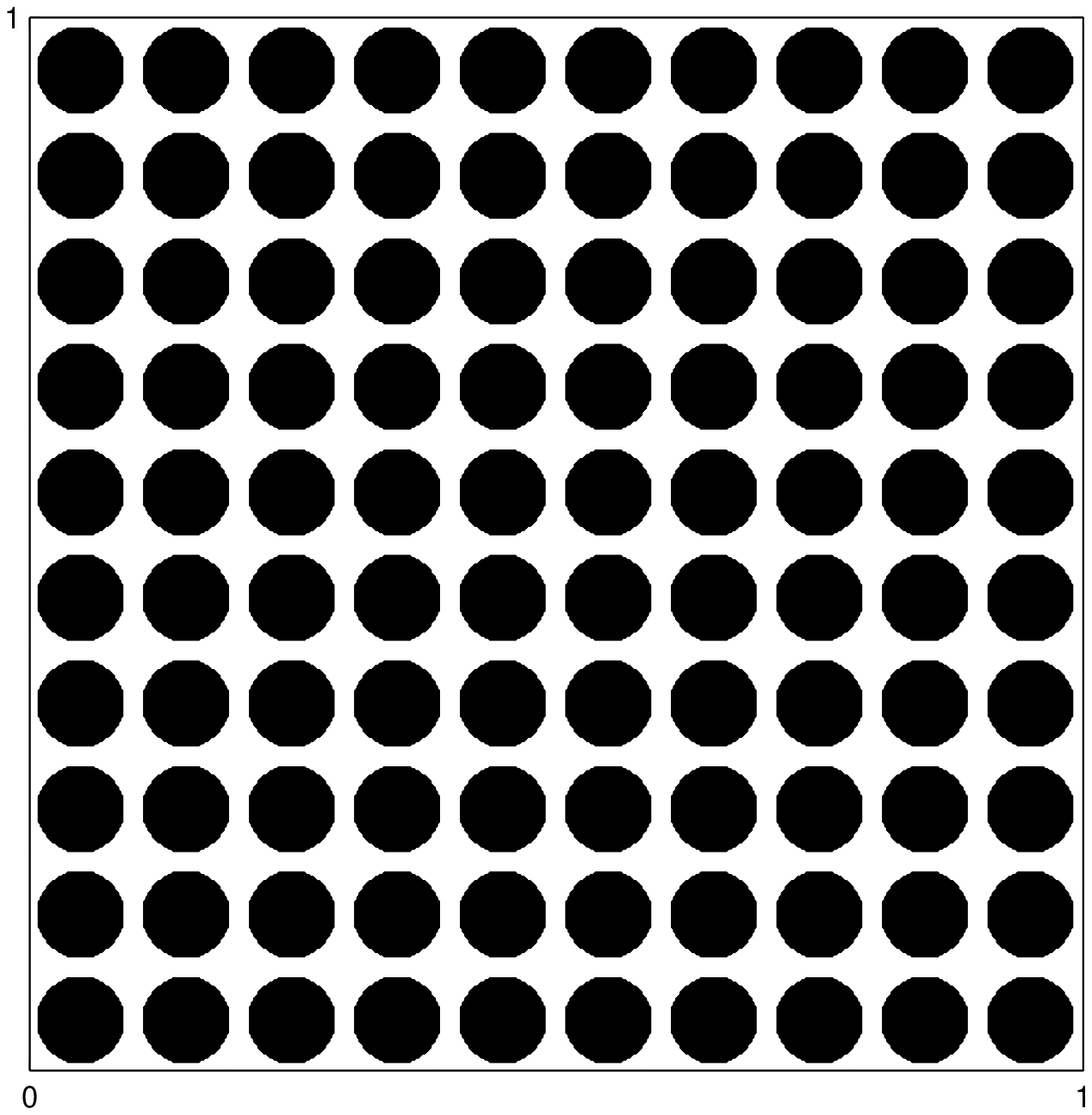}}
\subfigure[]{%
\label{fig+c} 
\includegraphics*[width=6.5cm, height=4.1
cm]{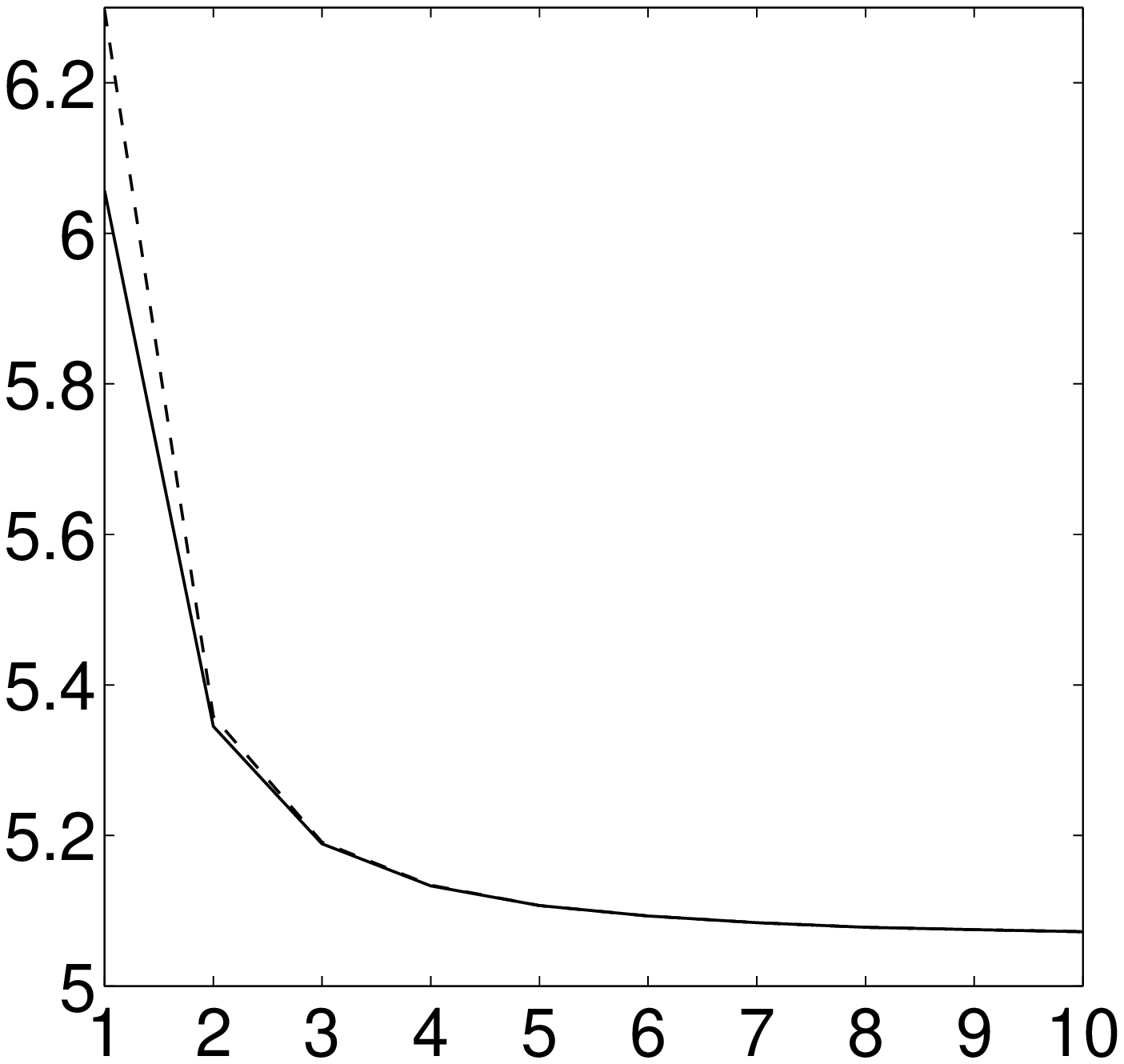}}
 \caption{$N=2$, $A=I_2$, $\nu\equiv 1$, $h\equiv 1$ and
$\mu(x)$ is $(1,1)$-periodic and takes two values,  $\mu\equiv-1$ on
$\Omega_-$ and $\mu\equiv 10$ on $\Omega^+$, where $\Omega^+$
consists on $k^2$ equally-spaced disks such that, on each period
cell $[0,1]^2$, $|\Omega^+\cap [0,1]^2|=1/2$, and
$\Omega^-=\R^2\backslash \Omega^+$. (a) A period cell with $k=2$ (b)
$k=10$; $\Omega^+$ is represented in black. (c) The values of
$\delta_1$ (continuous line) and $\delta_2$ (dashed line) in
function of $k$.}
 \label{fig}
\end{figure}

The proof relies on monotone methods of sub- and super-solutions.
For the existence result (i), we have computed a sub-solution of the
form $\kappa \phi$ with $\kappa > 0$. The optimal value of $\kappa$,
in the sense that it gives the highest value of $\delta_1$, is
$\kappa_0=\ds{-\lambda_1/(\overline{\nu}+\overline{\nu}\underline{\phi})}$.
We have numerically computed the values of $\delta_1$ and $\delta_2$
in several particular examples of sp-case (see Figure \ref{fig}).
The results illustrate the effect of environmental fragmentation on
the maximum sustainable yield, and show that the interval
$(\delta_1,\delta_2]$ on which we have no theoretical information
can be very narrow (see Figure \ref{fig}-(c)).

 Let us turn to the study of the
evolution equation (\ref{eq_evo2}). We assume that $\lambda_1<0$
 and $\epsilon$ is such that $\ds{\varepsilon_0:=2\frac{\varepsilon \overline{\nu}}{
\underline{\phi}}}<\ds{\frac{-\lambda_1}{2}}$;  we prove the
following theorem:

\vspace{1ex}

\begin{a2theorem} Let $u(t,x)$ be the solution of
{\rm(\ref{eq_evo2})} with initial data $u(0,x)=p(x)$ defined in
Remark \ref{rem1}. Then $u$ is non-increasing in $t$ and we have the
following asymptotic behaviour \par (i) if $\delta\leq \delta_1$,
$u(t,x)\to p_\delta(x)$ uniformly in $\Omega$ as $t\to +\infty$,
where $p_{\delta}$ is the unique positive maximal solution of
{\rm{(\ref{eq_sta})}}; and \par (ii) if $\delta
> \delta_2$, then $u(t,x)<  \varepsilon_0$
for $t$ large enough.
 \label{CRASth3}
\end{a2theorem}

\vspace{1ex}

\noindent In the above theorem, we assume that $u(0,x)=p(x).$ This
means that harvesting starts on a stabilized population governed by
the standard Fisher model without external forcing. \vspace{1ex}
\begin{a2remark}
{\rm These results are sharper than those which could be obtained by
a standard La Salle invariance principle, since we obtain here
discriminatory bounds on $\delta$, which determine the asymptotic
behaviour of the solutions.}
\end{a2remark}
\vspace{1ex}

\section{Time-periodic forcing}\label{sec3}

In this section we consider the general equation (\ref{eq_evo}) in
the {\bf bounded case}, with $\omega > 0$ defined as the frequency
of the forcing term. All the results are proved in \cite{cr}. Let us
introduce $T:=\omega^{-1}$. It is known that under the above
assumptions on $A(x)$, $\mathcal{A}u=-\nabla \cdot(A(x) \nabla u)$
is a sectorial operator with domain $\mathcal{D}(\mathcal{A})=\{u\in
H^2(\Omega), \mbox{ s. t. }
\partial_n u =0\ \mbox{ on } \partial \Omega \}$ (see e.g.
\cite{paz}). As a consequence, $-\mathcal{A}$ generates an analytic
semigroup $e^{-\mathcal{A}t}$ on $L^2(\Omega)$. Let $\{V^{2
r}\}_{\{r\geq 0\}}$ be the family of interpolation spaces generated
by the fractional powers of $\mathcal{A}$, where $V^{2
r}=\mathcal{D}(\mathcal{A}^{r})$ (see \cite{sell} for details). The
existence of a $T$-periodic solution of the equation (\ref{eq_evo})
can be reached by several procedures (e.g. averaging method
\cite{hal}). We present here a result on the existence of a
hyperbolic $T$-periodic solution, which is related to the robustness
of a hyperbolic stationary solution of the autonomous equation
(\ref{eq_evo2}), with $\delta h(x)= \int_{0}^1 f(s, x)ds$. More
precisely, \vspace{1ex}
\begin{a2theorem}\label{bogozthm}
Assume that equation {\rm(\ref{eq_evo2})} has a hyperbolic
stationary solution $q \in V^{2 r}$, $0\leq r \leq 1$. Then there
exists $\omega^{*} > 0$ such that for every $\omega \geq
\omega^{*}$, the problem {\rm(\ref{eq_evo})} possesses a hyperbolic
$T$-periodic solution $u_{\omega}(t)$ such that for any $t\in
[0,T]$, $u_{\omega}(t)$ lives in a $V^{2 r}$-neighborhood of $q$.
Furthermore if $\omega \rightarrow +\infty$, then $u_{\omega}(t)
\rightarrow q$ in $V^{2 r}$.
\end{a2theorem}
\vspace{1ex}

The proof uses similar arguments as the one of \cite{sell} Theorem
76.1, and is therefore based on a Lyapunov-Perron type argument; the
existence of such a hyperbolic periodic orbit is achieved via a
fixed point argument on the following operator:

\begin{equation*} \label{ope} \xi \rightarrow \mathcal{T}(f)\xi
:=\int_{-\infty}^t Q e^{-L(t-s)}(E(q,\xi)+f(ws,\cdot))ds...
\end{equation*}
\begin{equation*}
...-\int_t^{+\infty}P e^{-L(t-s)} (E(q,\xi)+f(ws,\cdot))ds,
\end{equation*}

where $L=\mathcal{A}-(\mu(x)-2\nu(x)q)\mathcal{I}$,
$E(q,\xi)=-2\nu(x)\xi^2$, $\xi$ belongs to a subset of
$L^{\infty}(\R, V^{2r})\cap C^0(\R,V^{2r})$,  and $P$ and $Q$ are
the associated projectors with the exponential dichotomy for
equation (\ref{eq_evo2}) related to the existence of a hyperbolic
stationary solution.

The main interest of Theorem \ref{bogozthm} is that it gives a
simple sufficient condition to ensure the existence of a
$T$-periodic solution of (\ref{eq_evo}) and that it allows to
localize in physical space where this solution can appear. Another
interesting aspect of this theorem is that it reduces the study of
existence and stability of a $T$-periodic solution of (\ref{eq_evo})
to that of the hyperbolic equilibria of the autonomous version
(\ref{eq_evo2}). For instance we get as an application of Theorems 2
and 4 for $\lambda_1 < 0$ and $f(wt,x)=\delta g(wt,x)$ with $\alpha<
\int_0^1 g(s,x)ds<\beta$ for all $x\in \Omega$, that if $\delta\leq
\delta_1$ and $\omega$ sufficiently large, then there exists a
stable non-trivial $T$-periodic solution $u_{\omega,\delta}$ of
$(\ref{eq_evo})$ in a neighborhood of a solution $p_{\delta}$ of
(\ref{eq_sta}).

\vspace{3ex}

{\bf Note added for this version on ArXiV:} The reference \cite{rc}
cited here, corresponds the article of the authors entitled ``On
Population Resilience to External Perturbations" which has been
published in {\it SIAM J. Appl. math (SIAP)}, {\bf 68} (1), (2007)
133–-153. We kept here the old reference \cite{rc} as in the
original article published in {\it C. R. Acad. Sci. Paris}, Ser. I,
343: 307-310; before the SIAP article. The reference \cite{cr} below
is still under preparation.


\begin{thebibliography}{99}


\footnotesize{
\bibitem{bhr1} H. Berestycki, F. Hamel, L. Roques. Analysis of the
periodically fragmented environment model~: I - Species persistence.
J. Math. Biol. 51 (1) (2005) 75-113.


\bibitem{cr}  M. Chekroun, L. Roques, Spatialized harvesting models. The influence of seasonal variations, in preparation.

\bibitem{fi} R.A. Fisher, The advance of advantageous genes, Ann.
Eugenics 7 (1937) 335-369.

\bibitem{hal} J.K. Hale, S.M. Verduyn Lunel, Averaging in Infinite
Dimensions, Journal of Int. Eq.  and Appl. 2 (4) (1990) 463-494.

\bibitem{paz} A. Pazy,
Semigroup of Linear Operators and Applications to Partial
Differential Equations, Springer-Verlag, 1983.


\bibitem{rc} L. Roques, M. Chekroun, Harvesting models in heterogeneous environments, Preprint (2006).

\bibitem{sell} G.R. Sell, Y. You, Dynamics of Evolutionary Equations, Springer-Verlag, 2002.

\bibitem{skt} N. Shigesada, K. Kawasaki, E. Teramoto, Traveling
periodic waves in heterogeneous environments, Theor. Population
Biol. 30 (1986) 143-160. }



\end{thebibliography}
\end{document}